# Connected Dominating Sets in Graphs with Stability Number Three

Vladimir Bercov

Department of Mathematics

Borough of Manhattan Community College

Abstract

In the special case of graphs $G$ of independence number $\alpha(G) = 3$ without induced chordless cycles $C_7$ it is shown that exists connected dominating set $D$ of vertices with number of vertices $n(D) \leq 4$. Using the concept of connected dominating sets, we defined a new invariant $h(G)$ that does not exceed the number of Hadwiger. For the considered graphs it is shown that $h(G) \geq n(G)/4$.

All graphs considered in this note are undirected, simple and finite. For a graph $G$, let $V(G)$ and $E(G)$ denote the vertex set and edge set of $G$. We will write $v \sim u$ ($v \nsim u$) when vertices $v$ and $u$ are (are not) adjacent. Further, let $n(G)$ and $\alpha(G)$ denote the number of vertices and the independence number of $G$ respectively. For $X \subseteq V(G)$, we denote by $G[X]$ the subgraph of $G$ induced by $X$, further $G - X = G[V(G) - X]$. The neighbourhood of vertex $v \in V(G)$, denoted $N(v)$, is a set of all vertices adjacent to $v$. The closed neighbourhood of $v$ is $N[v] = N(v) \cup \{v\}$. A simplicial vertex of a graph $G$ is a vertex $v$ for which $G[N(v)]$ is complete. $C_n$ is a chordless cycle with $n$ vertices: $C_n = v_1 v_2 \ldots v_n$. Connected dominating set in a graph $G = (V, E)$ is a subset $D \subseteq V$ such that every vertex $v \in V$ is in $D$ or $v$ is adjacent to some $u \in D$, and $G[D]$ is connected in all connected components of $G$. Using the concept of dominating sets, we can introduce a new invariant of a graph (see [1]). Number $h(G)$ can be defined as a maximum length of sequence $v_1, v_2, H_3, H_4, \ldots$, where $v_1$ and $v_2$ are adjacent vertices in graph $G$, $H_3$ is connected set of $V(G - \{v_1, v_2\})$ joined to $v_1$ and $v_2$; $H_4$ is connected set in $G - \{v_1, v_2\} - H_3$ joined to $v_1, v_2$ and to all vertices $H_3$, and so on. It is clear that $\omega(G) \leq h(G) \leq \eta(G)$, where $\omega(G)$ is a clique number and $\eta(G)$ is number of Hadwiger of a graph $G$, and for any connected dominating set $D$ in graph $G$: $h(G - D) \leq h(G) - 1$.

In [2] Plummer, Stiebitz and Toft proved : if $\alpha(G) \leq 2$ and G does not contain an induced $C_5$, then either $G$ contains complete $K_{n/2}$ or $G$ contains a dominating edge. We slightly modify this theorem.

**Theorem 1.** If graph $G$ is connected, $\alpha(G) = 2$ and $G$ does not contain an induced $C_5$ then for any non-simplicial vertex $v$ exists dominating edge $vu$.

*Proof.* Let $v \sim v_1$, $v \sim v_2$ and $v_1 \nsim v_2$; $G' = G - N[v]$. If $V(G') = \emptyset$ then $v$ is dominating vertex.

If $V(G') \neq \emptyset$ then denote by $V_1$ the set of neighbors of $v_1$ in $G'$ which are not neighbors of $v_2$, by $V_2$ the set of neighbors of $v_2$ in $G'$ which are not neighbors of $v_1$, by $V_{12}$ the set of vertices in $G'$ which are adjacent in $G'$ to $v_1$ and to $v_2$. Since $\alpha(G) = 2$, $G'$ is complete. If $V_1 \neq \emptyset$ and $V_2 \neq \emptyset$, then exists $C_5 = vv_1abv_2$, where $a \in V_1$, $b \in V_2$. If $V_1 = \emptyset$ (or $V_2 = \emptyset$) then $vv_2$ (or $vv_1$) is dominating edge.

**Corollary.** Let $G$ be a graph with $\alpha(G) = 2$. If $G$ does not contain an induced $C_5$ then $h(G) \geq n(G)/2$.

A graph is claw-free if no vertex has three pairwise nonadjacent neighbors.

**Theorem 2.** If graph $G$ is connected, claw-free, $\alpha(G) = 3$ and $G$ does not contain an induced $C_7$, then for any non-simplicial vertex $v$ exists connected dominating set $D$, such that $v \in D$ and $n(D) \leq 4$.

*Proof.* Let $v \sim v_1$, $v \sim v_2$ and $v_1 \nsim v_2$; $G' = G - N[v]$. Denote by $V_1$ the set of neighbors of $v_1$ in $G'$ which are not neighbors of $v_2$, by $V_2$ the set of neighbors of $v_2$ in $G'$ which are not neighbors of $v_1$, by $V_{12}$ the set of vertices in $G'$ which are adjacent in $G'$ to $v_1$ and to $v_2$, and by $V'_{12} = V(G') - (V_1 \cup V_2 \cup V_{12})$. Since $\alpha(G) = 3$, $G[V'_{12}]$ is complete or $V'_{12} = \emptyset$.

If $V'_{12} = \emptyset$ then $D = \{v, v_1, v_2\}$. If $V'_{12} \neq \emptyset$ and $V_1 \cup V_2 \cup V_{12} = \emptyset$ then $D = \{v, a, b\}$, where $a \sim v$, $b \sim a$ and $b \in V'_{12}$. Let $V'_{12} \neq \emptyset$ and $V_1 \cup V_2 \cup V_{12} \neq \emptyset$. In this case induced subgraphs $G[V_1 \cup V_{12}]$ and $G[V_2 \cup V_{12}]$ are complete (one of them can be null graph). Since $G$ is claw-free, all vertices $V_{12}$ do not have neighbors in $V'_{12}$. Since $\alpha(G) = 3$, $\alpha(G') \leq 2$.

Case 1. $G'$ is connected.

(I) $G'$ contains $C_5 = u_1u_2u_3u_4u_5$.
Since induced subgraphs $G[V_1 \cup V_{12}]$, $G[V_2 \cup V_{12}]$ and $G[V'_{12}]$ contain at most two consecutive vertices of $C_5$, there are three possible cases:

(a) $u_1 \in V_1$, $u_2 \in V_{12}$, $u_3 \in V_2$, $\{u_4, u_5\} \subseteq V'_{12}$.
In this case exists $C_7 = vv_2u_3u_4u_5u_1v_1$.

(b) $\{u_1, u_2\} \subseteq V_1$, $u_3 \in V_2$, $\{u_4, u_5\} \subseteq V'_{12}$.
In this case exists claw with the set of vertices $\{u_3, u_2, u_4, v_2\}$.

(c) $\{u_1, u_2\} \subseteq V_1, \{u_3, u_4\} \subseteq V_2, u_5 \in V'_{12}$.

Since $\alpha(G) = 3$, any pair of nonadjacent vertices from $V_1$ and $V_2$ is adjacent to all vertices from $V'_{12}$. If $u_4$ is connected to all vertices $V'_{12}$ then $D = \{v, v_1, v_2, u_4\}$. If $u \in V'_{12}$ and $u_4 \nsim u$, then $u_2 \sim u$, and in this case exists $C_7 = vv_2u_4u_5uu_2v_1$.

( II ) $G'$ does not contain $C_5$.

If $G'$ is complete then $D = \{v, a, b\}$, where $a \sim v$, $b \sim a$, $b \in V(G')$. Otherwise, by Th.1, in $G'$ for any non-simplicial vertex $b \in V(G')$ exists dominating edge $bc$. If $b$ is adjacent to $a \in N(v)$, then $D = \{v, a, b, c\}$.

Now let all vertices of $G'$ which are connected to $N(v)$ are simplicial.
If $V_{12} \neq \emptyset$ then $G[V_1 \cup V_2]$ and $G[V_1 \cup V_2 \cup V_{12}]$ are complete and, since $G'$ is connected, exists edge $ab$, where $b \in V'_{12}$ and $a \sim v_1$ (or $a \sim v_2$). In this case $D = \{v, v_1, a, b\}$ (or $D = \{v, v_2, a, b\}$).

Let $V_{12} = \emptyset$. Since all vertices of the sets $V_1$ and $V_2$ are simplicial, $G[V_1 \cup V_2]$ is complete (and we have the same dominating set as in previous case) or $V_1$ does not have neighbors in $V_2$. In the last case any two vertices $u_1 \in V_1$ and $u_2 \in V_2$ are adjacent to all vertices $V'_{12}$. If one of these two vertices is adjacent to all $V'_{12}$ then $D = \{v, v_1, v_2, u_1\}$ or $D = \{v, v_1, v_2, u_2\}$. Otherwise, exist $u_3, u_4 \in V'_{12}$ ($u_3$ and $u_4$ are different) such that $u_1 \nsim u_3$ and $u_2 \nsim u_4$. In this case exists $C_7 = vv_1u_1u_4u_3u_2v_2$.

Case 2. $G'$ is disconnected.

In this case $G'$ is disjoint union of two complete subgraphs $G'_1$ and $G'_2$ where $G'_1 = G[V_1 \cup V_2 \cup V_{12}]$, $G'_2 = G[V'_{12}]$. Since $G$ is claw-free, there are no vertices in $N(v)$ connected to both subgraphs. Denote by $U_1 \subseteq N(v)$ the subsetset of neighbors of $G'_1$ and by $U_2 \subseteq N(v)$ the subsetset of neighbors of $G'_2$. If $a \in U_1$, $b \in U_2$, and $a \nsim b$, then $a$ is connected to all $V(G'_1)$ or $b$ is connected to all $V(G'_2)$. There is no vertex $u_2 \in U_2$ which is connected to $v_1$ and to $v_2$. If $u_2 \nsim v_2$ and $v_2$ is adjacent to all vertices $V(G'_1)$ then $D = \{v, v_2, u_2, b\}$, where $b \in V'_{12}$ and $u_2 \sim b$. If $u_2 \nsim v_2$ and $u_2$ is adjacent to all vertices $V(G'_2)$ then $D = \{v, v_2, u_2, a\}$ where $a \in V_1 \cup V_2 \cup V_{12}$ and $v_2 \sim a$.

**Corollary.** Let $G$ be a graph with $\alpha(G) = 3$. If $G$ does not contain an induced $C_7$, then $h(G) \geq n(G)/4$.

*Proof.* We proceed by induction on $n = n(G)$. For $n \leq 4$, the result is clear. Suppose $n \geq 5$ and suppose the result is true for all graphs with fewer than $n$ vertices and let $G$ be a graph with $n$ vertices. If $G$ contains a claw, then the set $D$ of vertices of this claw is dominating in $G$, if not, by Th.2, we can build a dominating set $D$ with $n(D) \leq 4$. In both cases

$$h(G) \geq h(G - D) + 1 \geq \frac{n(G-D)}{4} + 1 \geq \frac{n(G)-4}{4} + 1 = n(G)/4.$$